\documentclass[12pt]{article}%
\usepackage{amsfonts}
\usepackage{sw20bams}
\usepackage{amsmath}
\usepackage{amssymb}
\usepackage{graphicx}%
\providecommand{\U}[1]{\protect\rule{.1in}{.1in}}
\begin{document}

\title{Stochastic Reservoir Calculations}
\author{Steven Finch}
\date{January 8, 2023}
\maketitle

\begin{abstract}
Prabhu (1958) obtained the stationary distribution of storage level $Z_{t}$ in
a reservoir of finite volume $v$, given an inflow $X_{t}$ and an outflow
$Y_{t}$. \ Time $t$ is assumed to be discrete, $X_{t}\sim$%
\ $\operatorname*{Gamma}(p,\mu)$ are independent and $p$ is a positive
integer. \ The mean inflow is $p/\mu$; the target outflow is $m$ (constant).
\ We attempt to clarify intricate details, often omitted in the literature, by
working through several examples. \ Of special interest are the probabilities
of depletion ($Z_{t}=0$) and spillage ($Z_{t}=v$). \ For prescribed
$\{v,p,\mu\}$, what value of $m$ minimizes both of these?

\end{abstract}

\footnotetext{Copyright \copyright \ 2023 by Steven R. Finch. All rights
reserved.}Let $v>0$, \thinspace$p$ be a positive integer, $\mu>0$ and $m>0$.
\ At each time $t=1,2,3,\ldots$, a reservoir of volume $v$ absorbs an inflow
$X_{t}\sim$\ $\operatorname*{Gamma}(p,\mu)$ and simultaneously releases an
outflow $0\leq Y_{t}\leq m$, depending on availablity. More precisely,
\[
Y_{t}=\min\left\{  X_{t}+Z_{t},m\right\}
\]
where $0\leq Z_{t}\leq v$ is the storage level. \ Independence across time is
assumed. \ Our interest is in the probability density function of $Z_{t}$ in
the limit as $t\rightarrow\infty$. \ We need not explicitly refer to $Y_{t}$
again, as $Z_{t+1}$ can be defined recursively without it:%
\[%
\begin{array}
[c]{ccc}%
Z_{t+1}=\max\left\{  0,\min\left\{  X_{t}+Z_{t}-m,v\right\}  \right\}  , &  &
Z_{1}=v/2\text{.}%
\end{array}
\]
Let $n=\left\lfloor v/m\right\rfloor $ and $\delta=v-m\,n$. \ In words,
$\delta$ is $0$ if and only if $v$ is an integer multiple of $m$, and $\delta$
is otherwise $>0$. \ Let%
\[
\lambda=(-1)^{p-1}\mu^{p}\exp(-\mu\,m).
\]
The graph of the PDF for $Z_{t}$ is piecewise smooth and contains at most
$n+1$ arcs, as well as point masses at $z=0$ and $z=v$. \ The arcs are
identified by $j=0,1,2,\ldots,n$ from left to right, and correspond to open
subintervals%
\[
\max\left\{  (j-1)m+\delta,0\right\}  <z<\min\left\{  j\,m+\delta,v\right\}
\]
of $0<z<v$. \ Prabhu \cite{Pra1-que4} impressively obtained the cumulative
distribution function%
\[
F(z)=1-\exp\left[  \mu(v-z)\right]
%TCIMACRO{\dsum \limits_{r=0}^{p-1}}%
%BeginExpansion
{\displaystyle\sum\limits_{r=0}^{p-1}}
%EndExpansion
\alpha_{r}%
%TCIMACRO{\dsum \limits_{q=0}^{n-j}}%
%BeginExpansion
{\displaystyle\sum\limits_{q=0}^{n-j}}
%EndExpansion
(-\lambda)^{q}\frac{(v-q\,m-z)^{q\,p+r}}{(q\,p+r)!}%
\]
that shall occupy us for the remainder of this paper. \ The $\alpha_{r}$
coefficients are found by solving a system of $p$ linear equations with
coefficients $d_{rs}$ for $r,s=0,1,\ldots,p-1$. \ These will be defined shortly.

Special considerations apply to endpoints. \ Let $\kappa=n-1$ if $\delta=0$
and $\kappa=n$ otherwise. \ The depletion probability, i.e., odds for the
reservoir to be dry, is%
\[
F(0)=1-\exp(\mu\,v)%
%TCIMACRO{\dsum \limits_{r=0}^{p-1}}%
%BeginExpansion
{\displaystyle\sum\limits_{r=0}^{p-1}}
%EndExpansion
\alpha_{r}%
%TCIMACRO{\dsum \limits_{q=0}^{\kappa}}%
%BeginExpansion
{\displaystyle\sum\limits_{q=0}^{\kappa}}
%EndExpansion
(-\lambda)^{q}\frac{(v-q\,m)^{q\,p+r}}{(q\,p+r)!}.
\]
In contrast, the spillage probability, i.e., odds for the reservoir to be
full, is just
\[
1-F(v)=\alpha_{0}.
\]
Minimizing the chance of both zero supply (harmful)\ and oversupply (wasteful)
is clearly important. \ Other quantities of interest include the total
deficit, i.e., unsatisfied demand, over a specified time duration; and total
surplus, i.e., unwanted supply (because $v<\infty$) that necessarily leaks
into the environment.

For $r=0,1,\ldots,p-1$, the linear system%
\[
\alpha_{r}-\lambda%
%TCIMACRO{\dsum \limits_{s=0}^{p-1}}%
%BeginExpansion
{\displaystyle\sum\limits_{s=0}^{p-1}}
%EndExpansion
d_{rs}\,\alpha_{s}=(-\mu)^{r}\exp\left[  -\mu(v+m)\right]
%TCIMACRO{\dsum \limits_{s=0}^{p-r-1}}%
%BeginExpansion
{\displaystyle\sum\limits_{s=0}^{p-r-1}}
%EndExpansion
\frac{\left[  \mu(v+m)\right]  ^{s}}{s!}%
\]
requires solution, where%
\[
d_{rs}=(-1)^{p+r-1}%
%TCIMACRO{\dsum \limits_{q=0}^{n}}%
%BeginExpansion
{\displaystyle\sum\limits_{q=0}^{n}}
%EndExpansion
(-\lambda)^{q}%
%TCIMACRO{\dint \limits_{q\,m}^{v}}%
%BeginExpansion
{\displaystyle\int\limits_{q\,m}^{v}}
%EndExpansion
\,\frac{(t-q\,m)^{q\,p+s}(t+m)^{p-r-1}}{(q\,p+s)!(p-r-1)!}\,dt.
\]
The integral can be easily expressed in closed-form.

Prabhu's CDF\ formula, given gamma-distributed inflow, extends a PDF\ formula
discovered earlier by Moran \cite{Mor1-que4}, given exponentially distributed
inflow ($p=1$). \ We have not studied \cite{Mor1-que4} in depth. \ More
discussion of \cite{Pra1-que4} appears in \cite{Mor2-que4, Pra2-que4, LP-que4,
Lyd-que4}. \ The treatment in \cite{Hm1-que4, Hm2-que4} is, however, most
pragmatic and useful for our purposes. \ 

Henceforth we fix $v=1$ and explore results for selected $\{p,\mu,m\}$. \ It
is surprising, more than fifty years after the publication of Prabhu's work,
that greater attention has not been paid to this research \cite{Hm1-que4}.
\ We can only imagine that intricate details, often lost in theoretical
summaries, have conspired to prevent greater understanding and widespread
recognition. \ Our hope is that working through a few examples will help to
improve matters.

\section{$\{p,\mu,m\}=\left\{  1,2,\frac{1}{2}\right\}  $}

The mean inflow is $p/\mu=1/2$ and the target outflow is $m=1/2$. \ Clearly
$n=\left\lfloor 1/m\right\rfloor =2$ and $\delta=1-m\,n=0$, i.e., there is no
offset. \ The arcs $j=0,1,2$ correspond to intervals%
\[%
\begin{array}
[c]{ccccc}%
0<z<0, &  & 0<z<1/2, &  & 1/2<z<1
\end{array}
\]
and thus $j=0$ can be ignored (being empty). \ Prabhu's formula gives $F(z)$
as%
\[
1-\frac{1}{2}\exp\left[  \mu(1-z)\right]  \left[  2-\left(  1-2z\right)
\lambda\right]  \alpha_{0}%
\]
for $j=1$ and
\[
1-\exp\left[  \mu(1-z)\right]  \alpha_{0}%
\]
for $j=2$. \ The linear equation
\[
(1-\lambda\,d_{00})\alpha_{0}=\exp\left(  -\frac{3}{2}\mu\right)
\]
coupled with%
\[
d_{00}=1-\frac{1}{8}\lambda
\]
and $\lambda=2e^{-1}$ give
\[
\alpha_{0}=\frac{8}{8-8\lambda+\lambda^{2}}\exp\left(  -\frac{3}{2}\mu\right)
=0.15000227...
\]
as the spillage probability. \ Because $\kappa=n-1=1$,
\[
F(0)=1-\frac{1}{2}\exp(\mu)\left(  2-\lambda\right)  \alpha_{0}=0.29937324...
\]
is the depletion probability. \ One may have expected these two probabilities
to be almost equal (since $1/\mu=1/2=m$ and by a certain symmetry), but this
is not true. \ The derivative $f(z)$ of $F(z)$ is plotted in Figure 1.

\section{$\{p,\mu,m\}=\left\{  1,2,\frac{1}{3}\right\}  $}

The mean inflow is $p/\mu=1/2$ and the target outflow is $m=1/3$. \ Clearly
$n=\left\lfloor 1/m\right\rfloor =3$ and $\delta=1-m\,n=0$, i.e., there is no
offset. \ The arcs $j=0,1,2,3$ correspond to intervals%
\[%
\begin{array}
[c]{ccccccc}%
0<z<0, &  & 0<z<1/3, &  & 1/3<z<2/3, &  & 2/3<z<1
\end{array}
\]
and thus $j=0$ can be ignored (being empty). \ Prabhu's formula gives $F(z)$
as%
\[
1-\frac{1}{18}\exp\left[  \mu(1-z)\right]  \left[  18-6\left(  2-3z\right)
\lambda+(1-3z)^{2}\lambda^{2}\right]  \alpha_{0}%
\]
for $j=1$,
\[
1-\frac{1}{3}\exp\left[  \mu(1-z)\right]  \left[  3-\left(  2-3z\right)
\lambda\right]  \alpha_{0}%
\]
for $j=2$ and%
\[
1-\exp\left[  \mu(1-z)\right]  \alpha_{0}%
\]
for $j=3$. \ The linear equation
\[
(1-\lambda\,d_{00})\alpha_{0}=\exp\left(  -\frac{4}{3}\mu\right)
\]
coupled with%
\[
d_{00}=1-\frac{2}{9}\lambda+\frac{1}{162}\lambda^{2}%
\]
and $\lambda=2e^{-2/3}$ give
\[
\alpha_{0}=\frac{162}{162-162\lambda+36\lambda^{2}-\lambda^{3}}\exp\left(
-\frac{4}{3}\mu\right)  =0.34604845...
\]
as the spillage probability. \ Because $\kappa=n-1=1$,%
\[
F(0)=1-\frac{1}{18}\exp(\mu)\left(  18-12\lambda+\lambda^{2}\right)
\alpha_{0}=0.04363903...
\]
is the depletion probability. \ While $\alpha_{0}<F(0)$ in Section 1, we have
$\alpha_{0}>F(0)$ here. \ This outcome suggests examining a value of $m$
between $1/3$ and $1/2$. \ The derivative $f(z)$ of $F(z)$ is plotted in
Figure 2.

\section{$\{p,\mu,m\}=\left\{  1,2,\frac{2}{5}\right\}  $}

The mean inflow is $p/\mu=1/2$ and the target outflow is $m=2/5$. \ Clearly
$n=\left\lfloor 1/m\right\rfloor =2$ and $\delta=1-m\,n=1/5$, i.e., the offset
is nonzero. \ The arcs $j=0,1,2$ correspond to intervals%
\[%
\begin{array}
[c]{ccccc}%
0<z<1/5, &  & 1/5<z<3/5, &  & 3/5<z<1;
\end{array}
\]
note that $j=0$ has length only $1/5$. \ Prabhu's formula gives $F(z)$ as%
\[
1-\frac{1}{50}\exp\left[  \mu(1-z)\right]  \left[  50-10\left(  3-5z\right)
\lambda+(1-5z)^{2}\lambda^{2}\right]  \alpha_{0}%
\]
for $j=0$,%
\[
1-\frac{1}{5}\exp\left[  \mu(1-z)\right]  \left[  5-\left(  3-5z\right)
\lambda\right]  \alpha_{0}%
\]
for $j=1$ and
\[
1-\exp\left[  \mu(1-z)\right]  \alpha_{0}%
\]
for $j=2$. \ The linear equation
\[
(1-\lambda\,d_{00})\alpha_{0}=\exp\left(  -\frac{7}{5}\mu\right)
\]
coupled with
\[
d_{00}=1-\frac{9}{50}\lambda+\frac{1}{750}\lambda^{2}%
\]
and $\lambda=2e^{-4/5}$ give
\[
\alpha_{0}=\frac{750}{750-750\lambda+135\lambda^{2}-\lambda^{3}}\exp\left(
-\frac{7}{5}\mu\right)  =0.24745701...
\]
as the spillage probability. \ Because $\kappa=n=2$,%
\[
F(0)=1-\frac{1}{50}\exp(\mu)\left(  50-30\lambda+\lambda^{2}\right)
\alpha_{0}=0.12789671...
\]
is the depletion probability. \ The values $\alpha_{0}$ and $F(0)$ are closer
than in the previous two sections; a choice of $m$ that is intermediate to
$2/5$ and $1/2$ should make these coincident. \ We estimate that $m=0.44276$
meets this objective (with $0.199$ as the common probability). \ On the other
hand, if our goal is to minimize the unweighted combination $\alpha_{0}+F(0)$,
then $m=0.38$ achieves the goal (with sum $0.372$). \ The derivative $f(z)$ of
$F(z)$ is plotted in Figure 3.

\section{$\{p,\mu,m\}=\left\{  2,4,\frac{1}{2}\right\}  $}

The mean inflow is $p/\mu=1/2$ and the target outflow is $m=1/2$. \ Clearly
$n=\left\lfloor 1/m\right\rfloor =2$ and $\delta=1-m\,n=0$, i.e., there is no
offset. \ The arcs $j=0,1,2$ correspond to intervals%
\[%
\begin{array}
[c]{ccccc}%
0<z<0, &  & 0<z<1/2, &  & 1/2<z<1
\end{array}
\]
and thus $j=0$ can be ignored (being empty). \ Prabhu's formula gives $F(z)$
as%
\[
1-\frac{1}{48}\exp\left[  \mu(1-z)\right]  \left\{  \left[  48-6\left(
1-2z\right)  ^{2}\lambda\right]  \alpha_{0}+\left[  48-48z-\left(
1-2z\right)  ^{3}\lambda\right]  \alpha_{1}\right\}
\]
for $j=1$ and%
\[
1-\exp\left[  \mu(1-z)\right]  \left\{  \alpha_{0}+(1-z)\alpha_{1}\right\}
\]
for $j=2$. \ The linear equations
\[
(1-\lambda\,d_{00})\alpha_{0}-\lambda\,d_{01}\alpha_{1}=\exp\left(  -\frac
{3}{2}\mu\right)  \left(  1+\frac{3}{2}\mu\right)  ,
\]%
\[
\lambda\,d_{10}\alpha_{0}-(1-\lambda\,d_{11})\alpha_{1}=\exp\left(  -\frac
{3}{2}\mu\right)  \mu
\]
coupled with%
\[%
\begin{array}
[c]{ccc}%
d_{00}=-1+\dfrac{11}{384}\lambda, &  & d_{01}=-\dfrac{7}{12}+\dfrac{7}%
{1920}\lambda,
\end{array}
\]%
\[%
\begin{array}
[c]{ccc}%
d_{10}=1-\dfrac{1}{48}\lambda, &  & d_{11}=\dfrac{1}{2}-\dfrac{1}{384}\lambda
\end{array}
\]
and $\lambda=-16e^{-2}$ give
\[
\alpha_{0}=\frac{1}{2}\frac{2+3\mu-2\lambda\,\mu\,d_{01}-\lambda\left(
2+3\mu\right)  d_{11}}{1-\lambda\left(  d_{00}+d_{11}\right)  +\lambda
^{2}\left(  d_{00}d_{11}-d_{01}d_{10}\right)  }\exp\left(  -\frac{3}{2}%
\mu\right)  ,
\]%
\[
\alpha_{1}=\frac{1}{2}\frac{-2\mu+2\lambda\,\mu\,d_{00}+\lambda\left(
2+3\mu\right)  d_{10}}{1-\lambda\left(  d_{00}+d_{11}\right)  +\lambda
^{2}\left(  d_{00}d_{11}-d_{01}d_{10}\right)  }\exp\left(  -\frac{3}{2}%
\mu\right)  ;
\]
the spillage probability is hence $\alpha_{0}=0.13554701...$. \ Because
$\kappa=n-1=1$,%
\[
F(0)=1-\frac{1}{48}\exp(\mu)\left[  \left(  48-6\lambda\right)  \alpha
_{0}+(48-\lambda)\alpha_{1}\right]  =0.22163253...
\]
is the depletion probability. \ The mode of $\operatorname*{Gamma}(2,\mu)$ is
$1/\mu>0$ whereas the mode of $\operatorname*{Gamma}(1,\mu)$ is $0$; a small
inflow is less likely for $p=2$ than for $p=1$, thus $F(0)$ is noticeably
smaller than in Section 1. \ The tail of $\operatorname*{Gamma}(2,\mu)$ is
fatter than the tail of $\operatorname*{Gamma}(1,\mu)$; a large inflow is more
likely for $p=2$ than for $p=1$, however $\alpha_{0}$ is paradoxically smaller
than in Section 1 (but only slightly). \ The derivative $f(z)$ of $F(z)$ is
plotted in Figure 4.

\section{Invariance}

One verification of Prabhu's formula is based on simulation (easily
programmed, since the recurrence for $Z_{t}$ is straightforward). \ Another
verification is more esoteric: to confirm that the formula is invariant under
the transformation%
\[
\left\{  v,\frac{p}{\mu},m\right\}  \longmapsto\left\{  \tilde{v},\frac
{p}{\tilde{\mu}},\tilde{m}\right\}  =\left\{  \frac{v}{m},\frac{p}{m\,\mu
},1\right\}
\]
in the sense that spillage \& depletion probabilities should remain constant
and storage level CDF arguments should simply scale by $m$. \ First,%
\[
\tilde{n}=\left\lfloor \frac{\tilde{v}}{\tilde{m}}\right\rfloor =\left\lfloor
\frac{v}{m}\right\rfloor =n,
\]%
\[
\tilde{\lambda}=(-1)^{p-1}\tilde{\mu}^{p}\exp[-\tilde{\mu}\,\tilde
{m}]=(-1)^{p-1}(m\,\mu)^{p}\exp[-m\,\mu\cdot1]=m^{p}\lambda
\]
and%
\begin{align*}
\tilde{d}_{rs}  &  =(-1)^{p+r-1}%
%TCIMACRO{\dsum \limits_{q=0}^{n}}%
%BeginExpansion
{\displaystyle\sum\limits_{q=0}^{n}}
%EndExpansion
(-\tilde{\lambda})^{q}%
%TCIMACRO{\dint \limits_{q\,\tilde{m}}^{\tilde{v}}}%
%BeginExpansion
{\displaystyle\int\limits_{q\,\tilde{m}}^{\tilde{v}}}
%EndExpansion
\,\frac{(t-q\,\tilde{m})^{q\,p+s}(t+\tilde{m})^{p-r-1}}{(q\,p+s)!(p-r-1)!}%
\,dt\\
&  =(-1)^{p+r-1}%
%TCIMACRO{\dsum \limits_{q=0}^{n}}%
%BeginExpansion
{\displaystyle\sum\limits_{q=0}^{n}}
%EndExpansion
m^{p\,q}(-\lambda)^{q}%
%TCIMACRO{\dint \limits_{q}^{v/m}}%
%BeginExpansion
{\displaystyle\int\limits_{q}^{v/m}}
%EndExpansion
\,\frac{(t-q)^{q\,p+s}(t+1)^{p-r-1}}{(q\,p+s)!(p-r-1)!}\,dt\\
&  =(-1)^{p+r-1}%
%TCIMACRO{\dsum \limits_{q=0}^{n}}%
%BeginExpansion
{\displaystyle\sum\limits_{q=0}^{n}}
%EndExpansion
m^{p\,q}(-\lambda)^{q}%
%TCIMACRO{\dint \limits_{q\,m}^{v}}%
%BeginExpansion
{\displaystyle\int\limits_{q\,m}^{v}}
%EndExpansion
\,\frac{(\frac{u}{m}-q)^{q\,p+s}(\frac{u}{m}+1)^{p-r-1}}{(q\,p+s)!(p-r-1)!}%
\,\frac{du}{m}%
\end{align*}
upon setting $u=m\,t$, $du=m\,dt$; thus%
\begin{align*}
\tilde{d}_{rs}  &  =(-1)^{p+r-1}%
%TCIMACRO{\dsum \limits_{q=0}^{n}}%
%BeginExpansion
{\displaystyle\sum\limits_{q=0}^{n}}
%EndExpansion
\frac{m^{p\,q}(-\lambda)^{q}}{m^{p\,q+s+p-r-1+1}}%
%TCIMACRO{\dint \limits_{q\,m}^{v}}%
%BeginExpansion
{\displaystyle\int\limits_{q\,m}^{v}}
%EndExpansion
\,\frac{(u-q\,m)^{q\,p+s}(u+m)^{p-r-1}}{(q\,p+s)!(p-r-1)!}\,du\\
&  =m^{-(p-r+s)}d_{rs}.
\end{align*}
Second,%
\[
\tilde{\alpha}_{r}-\tilde{\lambda}%
%TCIMACRO{\dsum \limits_{s=0}^{p-1}}%
%BeginExpansion
{\displaystyle\sum\limits_{s=0}^{p-1}}
%EndExpansion
\tilde{d}_{rs}\,\tilde{\alpha}_{s}=(-\tilde{\mu})^{r}\exp\left[  -\tilde{\mu
}(\tilde{v}+\tilde{m})\right]
%TCIMACRO{\dsum \limits_{s=0}^{p-r-1}}%
%BeginExpansion
{\displaystyle\sum\limits_{s=0}^{p-r-1}}
%EndExpansion
\frac{\left[  \tilde{\mu}(\tilde{v}+\tilde{m})\right]  ^{s}}{s!}%
\]
implies%
\[
\tilde{\alpha}_{r}-m^{p}\lambda%
%TCIMACRO{\dsum \limits_{s=0}^{p-1}}%
%BeginExpansion
{\displaystyle\sum\limits_{s=0}^{p-1}}
%EndExpansion
m^{-(p-r+s)}d_{rs}\,\tilde{\alpha}_{s}=m^{r}(-\mu)^{r}\exp\left[
-\mu(v+m)\right]
%TCIMACRO{\dsum \limits_{s=0}^{p-r-1}}%
%BeginExpansion
{\displaystyle\sum\limits_{s=0}^{p-r-1}}
%EndExpansion
\frac{\left[  \mu(v+m)\right]  ^{s}}{s!}%
\]
because $\tilde{\mu}(\tilde{v}+\tilde{m})=m\,\mu\left(  \frac{v}{m}+1\right)
=\mu(v+m)$; therefore%
\[
m^{-r}\tilde{\alpha}_{r}-\lambda%
%TCIMACRO{\dsum \limits_{s=0}^{p-1}}%
%BeginExpansion
{\displaystyle\sum\limits_{s=0}^{p-1}}
%EndExpansion
m^{-s}d_{rs}\,\tilde{\alpha}_{s}=(-\mu)^{r}\exp\left[  -\mu(v+m)\right]
%TCIMACRO{\dsum \limits_{s=0}^{p-r-1}}%
%BeginExpansion
{\displaystyle\sum\limits_{s=0}^{p-r-1}}
%EndExpansion
\frac{\left[  \mu(v+m)\right]  ^{s}}{s!}%
\]
which is immediately satisfied by $\tilde{\alpha}_{r}=m^{r}\alpha_{r}$. \ In
particular, $\tilde{\alpha}_{0}=\alpha_{0}$. \ Third,%
\[
\tilde{\delta}=\tilde{v}-\tilde{m}\,n=\frac{v-m\,n}{m}=\frac{\delta}{m}.
\]
Finally, given $j$,%
\begin{align*}
\tilde{F}(z)  &  =1-\exp\left[  \tilde{\mu}(\tilde{v}-z)\right]
%TCIMACRO{\dsum \limits_{r=0}^{p-1}}%
%BeginExpansion
{\displaystyle\sum\limits_{r=0}^{p-1}}
%EndExpansion
\tilde{\alpha}_{r}%
%TCIMACRO{\dsum \limits_{q=0}^{n-j}}%
%BeginExpansion
{\displaystyle\sum\limits_{q=0}^{n-j}}
%EndExpansion
(-\tilde{\lambda})^{q}\frac{(\tilde{v}-q\,\tilde{m}-z)^{q\,p+r}}{(q\,p+r)!}\\
&  =1-\exp\left[  (m\,\mu)\left(  \frac{v}{m}-z\right)  \right]
%TCIMACRO{\dsum \limits_{r=0}^{p-1}}%
%BeginExpansion
{\displaystyle\sum\limits_{r=0}^{p-1}}
%EndExpansion
m^{r}\alpha_{r}%
%TCIMACRO{\dsum \limits_{q=0}^{n-j}}%
%BeginExpansion
{\displaystyle\sum\limits_{q=0}^{n-j}}
%EndExpansion
m^{p\,q}(-\lambda)^{q}\frac{(\frac{v}{m}-q-z)^{q\,p+r}}{(q\,p+r)!}\\
&  =1-\exp\left[  \mu(v-m\,z)\right]
%TCIMACRO{\dsum \limits_{r=0}^{p-1}}%
%BeginExpansion
{\displaystyle\sum\limits_{r=0}^{p-1}}
%EndExpansion
m^{r}\alpha_{r}%
%TCIMACRO{\dsum \limits_{q=0}^{n-j}}%
%BeginExpansion
{\displaystyle\sum\limits_{q=0}^{n-j}}
%EndExpansion
\frac{m^{p\,q}(-\lambda)^{q}}{m^{p\,q+r}}\frac{(v-q\,m-m\,z)^{q\,p+r}%
}{(q\,p+r)!}\\
&  =F(m\,z)
\end{align*}
for $(j-1)\tilde{m}+\tilde{\delta}<z<j\,\tilde{m}+\tilde{\delta}$, i.e.,
$(j-1)m+\delta<m\,z<j\,m+\delta$. \ In the same way, $\tilde{F}(0)=F(0)$, with
the upper summation limit $n-j$ replaced by $\tilde{\kappa}=\kappa$.%
%TCIMACRO{\FRAME{ftbpFU}{4.3846in}{2.8496in}{0pt}{\Qcb{Plot of storage level
%density $w=f(z)$ for $\{p,\mu,m\}=\left\{  1,2,\frac{1}{2}\right\}  $.}}%
%{}{figur1.eps}{\special{ language "Scientific Word";  type "GRAPHIC";
%width 4.3846in;  height 2.8496in;  depth 0pt;  original-width 13.255in;
%original-height 8.5867in;  cropleft "0";  croptop "1";  cropright "1";
%cropbottom "0";  filename 'figur1.eps';file-properties "XNPEU";}}}%
%BeginExpansion
\begin{figure}
[ptb]
\begin{center}
\includegraphics[
height=2.8496in,
width=4.3846in
]%
{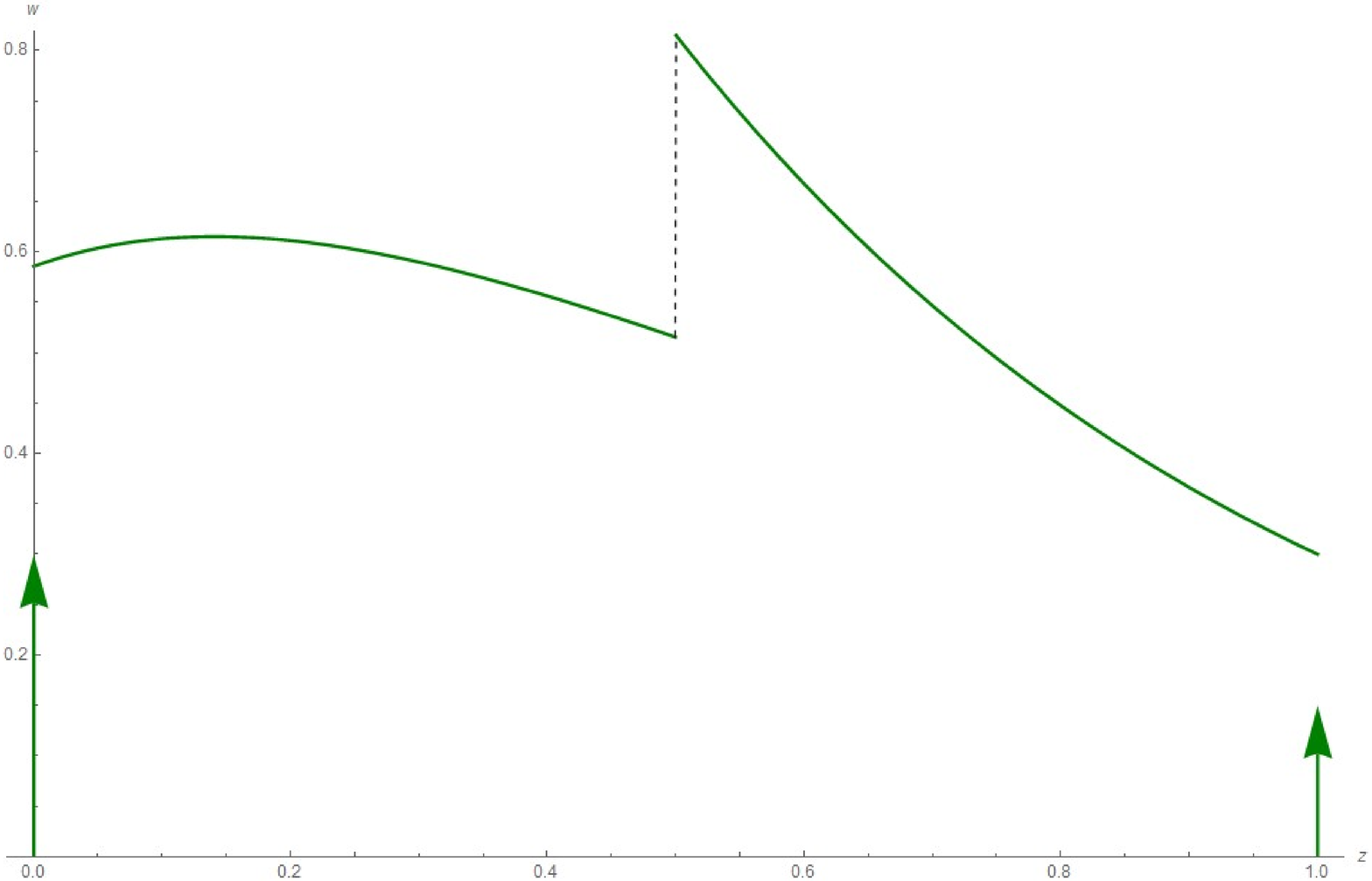}%
\caption{Plot of storage level density $w=f(z)$ for $\{p,\mu,m\}=\left\{
1,2,\frac{1}{2}\right\}  $.}%
\end{center}
\end{figure}
%EndExpansion%
%TCIMACRO{\FRAME{ftbpFU}{4.3759in}{2.8504in}{0pt}{\Qcb{Plot of storage level
%density $w=f(z)$ for $\{p,\mu,m\}=\left\{  1,2,\frac{1}{3}\right\}  $.}}%
%{}{figur2.eps}{\special{ language "Scientific Word";  type "GRAPHIC";
%width 4.3759in;  height 2.8504in;  depth 0pt;  original-width 12.7543in;
%original-height 8.2806in;  cropleft "0";  croptop "1";  cropright "1";
%cropbottom "0";  filename 'figur2.eps';file-properties "XNPEU";}}}%
%BeginExpansion
\begin{figure}
[ptb]
\begin{center}
\includegraphics[
height=2.8504in,
width=4.3759in
]%
{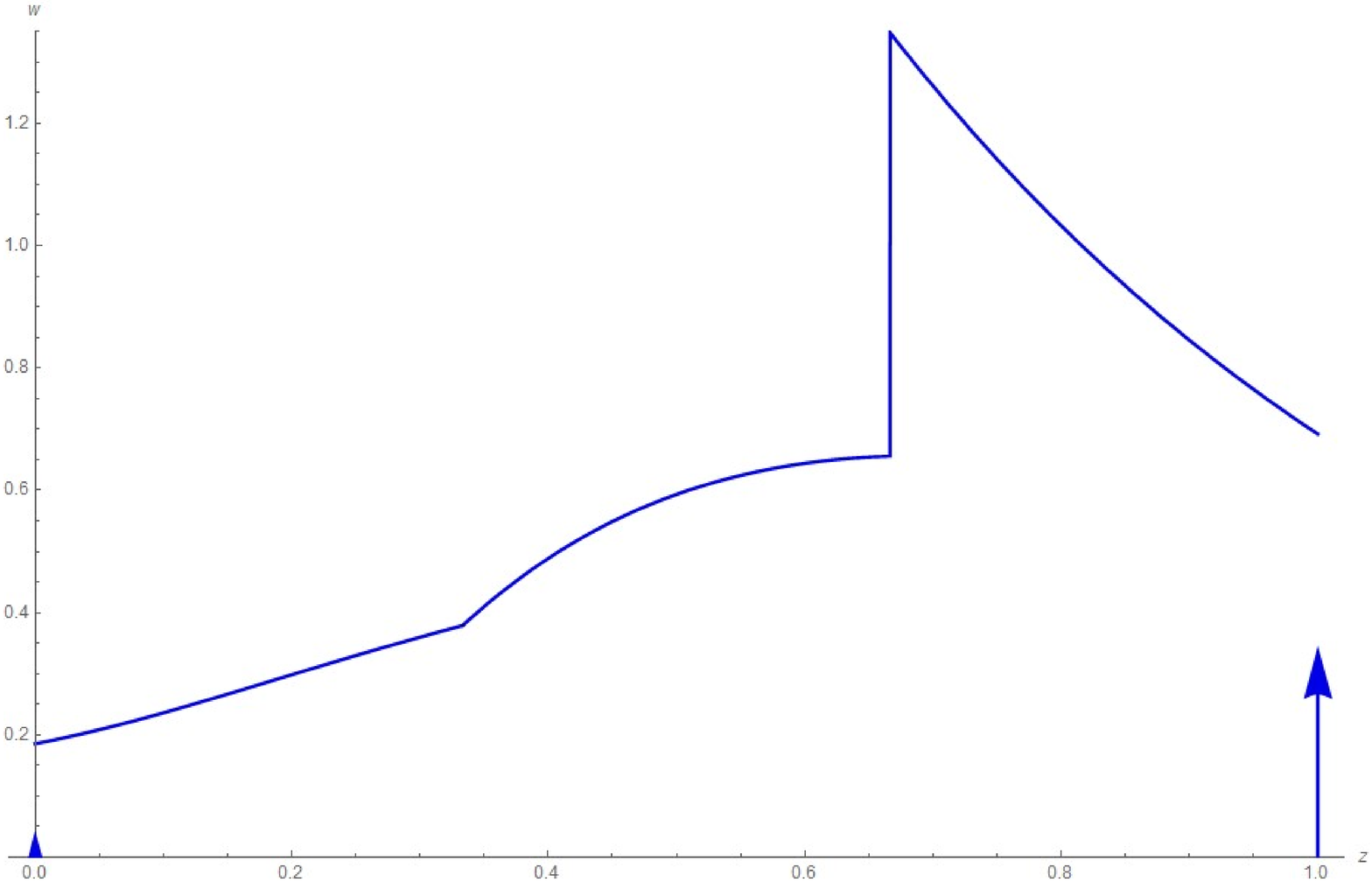}%
\caption{Plot of storage level density $w=f(z)$ for $\{p,\mu,m\}=\left\{
1,2,\frac{1}{3}\right\}  $.}%
\end{center}
\end{figure}
%EndExpansion%
%TCIMACRO{\FRAME{ftbpFU}{4.3803in}{2.8496in}{0pt}{\Qcb{Plot of storage level
%density $w=f(z)$ for $\{p,\mu,m\}=\left\{  1,2,\frac{2}{5}\right\}  $.}}%
%{}{figur3.eps}{\special{ language "Scientific Word";  type "GRAPHIC";
%width 4.3803in;  height 2.8496in;  depth 0pt;  original-width 12.8364in;
%original-height 8.323in;  cropleft "0";  croptop "1";  cropright "1";
%cropbottom "0";  filename 'figur3.eps';file-properties "XNPEU";}}}%
%BeginExpansion
\begin{figure}
[ptb]
\begin{center}
\includegraphics[
height=2.8496in,
width=4.3803in
]%
{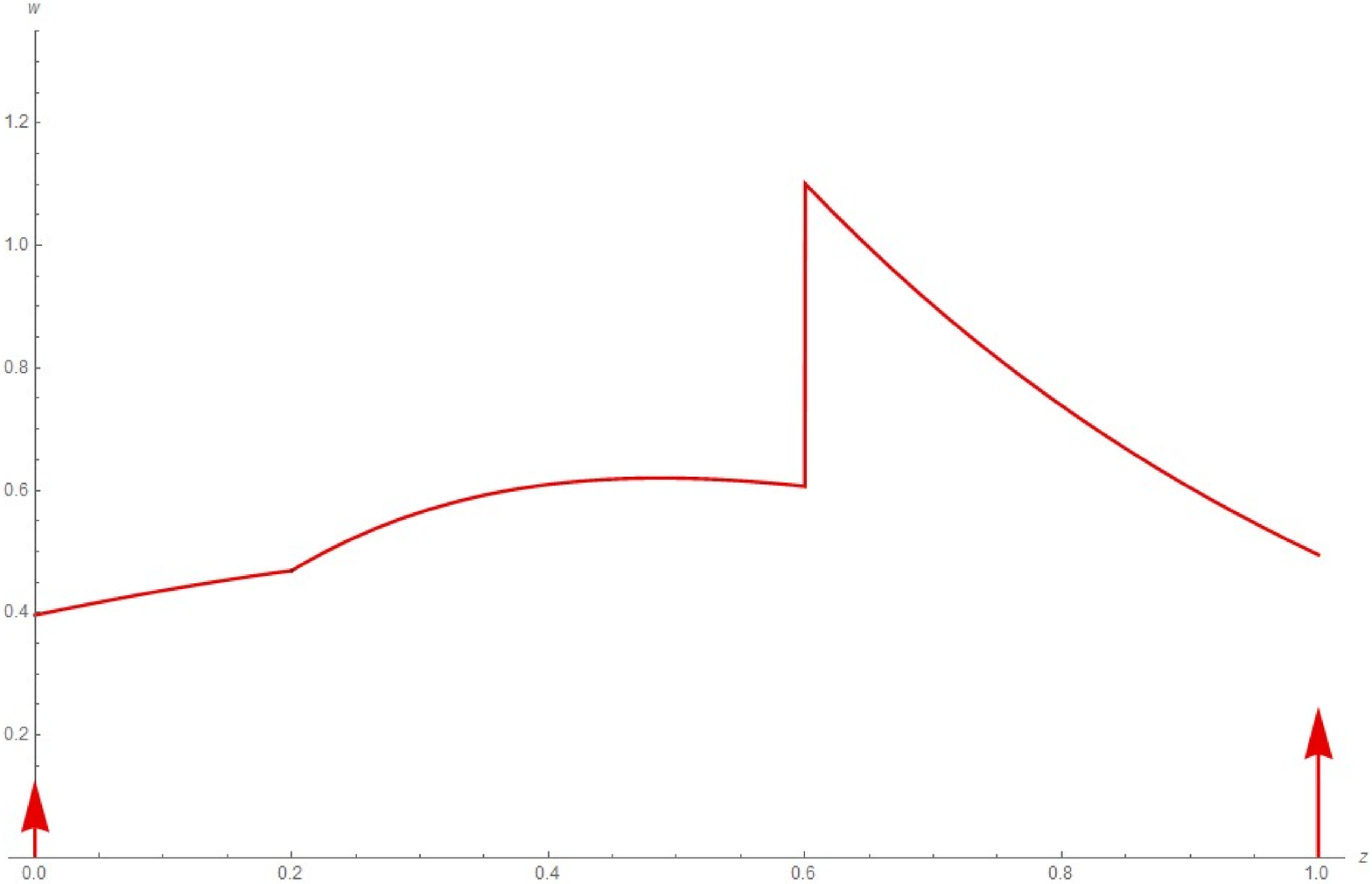}%
\caption{Plot of storage level density $w=f(z)$ for $\{p,\mu,m\}=\left\{
1,2,\frac{2}{5}\right\}  $.}%
\end{center}
\end{figure}
%EndExpansion%
%TCIMACRO{\FRAME{ftbpFU}{4.3803in}{2.8504in}{0pt}{\Qcb{Plot of storage level
%density $w=f(z)$ for $\{p,\mu,m\}=\left\{  2,4,\frac{1}{2}\right\}  $.}}%
%{}{figur4.eps}{\special{ language "Scientific Word";  type "GRAPHIC";
%width 4.3803in;  height 2.8504in;  depth 0pt;  original-width 12.8096in;
%original-height 8.3083in;  cropleft "0";  croptop "1";  cropright "1";
%cropbottom "0";  filename 'figur4.eps';file-properties "XNPEU";}}}%
%BeginExpansion
\begin{figure}
[ptb]
\begin{center}
\includegraphics[
height=2.8504in,
width=4.3803in
]%
{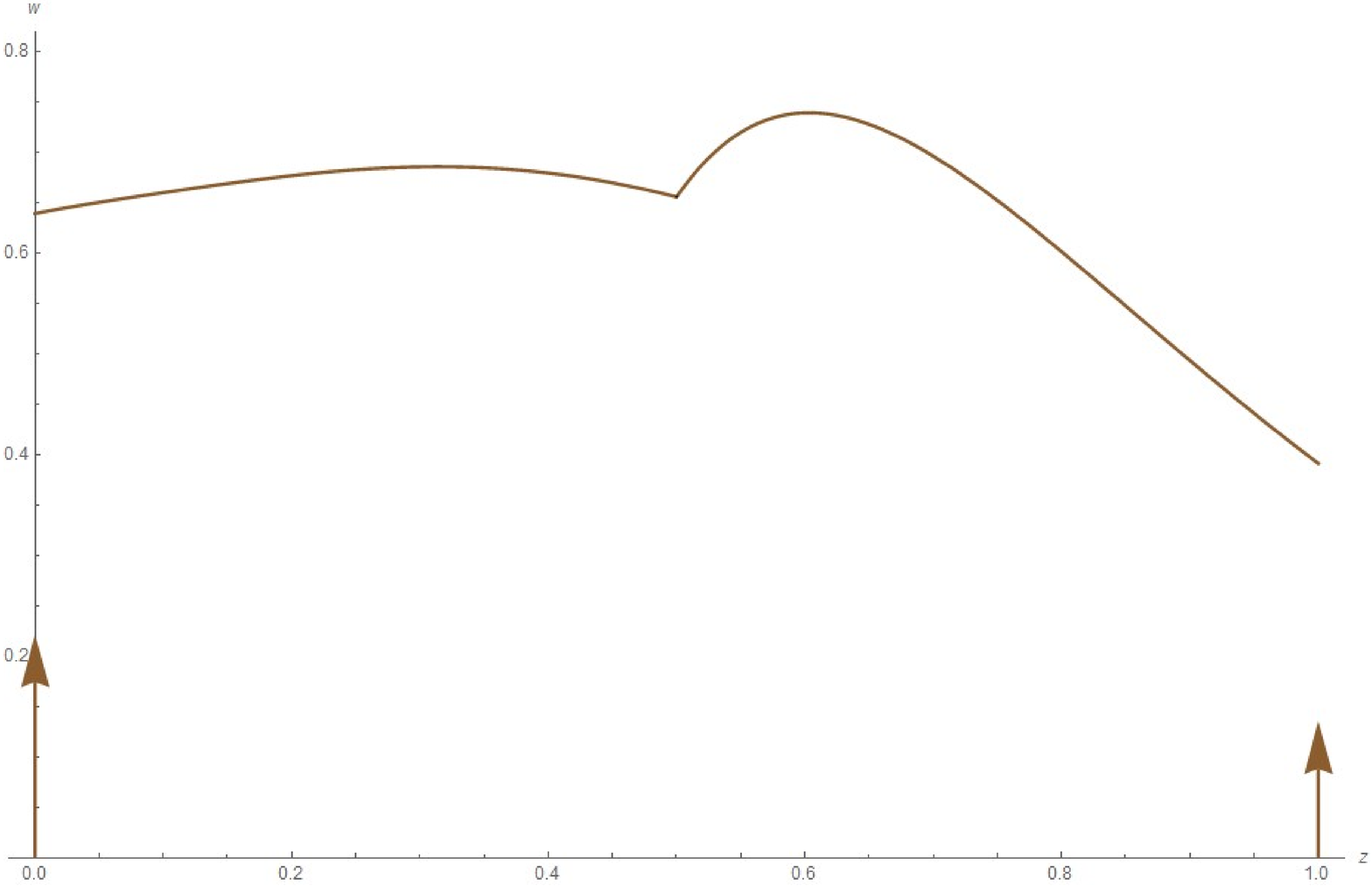}%
\caption{Plot of storage level density $w=f(z)$ for $\{p,\mu,m\}=\left\{
2,4,\frac{1}{2}\right\}  $.}%
\end{center}
\end{figure}
%EndExpansion

\section{Inquiry}

Moran \cite{Mor3-que4, Gani-que4} introduced a different model -- in
continuous time -- for an infinite volume reservoir. \ Let $X(t)\sim
\operatorname*{Gamma}(t,1/\rho)$ denote the total inflow over the interval
$(0,t]$, assumed to be a nonnegative stochastic process with stationary
independent increments, where $0<\rho<1$ is constant. \ In particular,
$\mathbb{E}(X(T))=\rho\,t$. \ Let the outflow be continuous and at unit rate
except when the reservoir is empty. \ We have%
\[
Z(t)=Z(0)+X(t)-t+%
%TCIMACRO{\dint \limits_{0}^{t}}%
%BeginExpansion
{\displaystyle\int\limits_{0}^{t}}
%EndExpansion
1_{\left\{  Z(\tau)=0\right\}  }\,d\tau
\]
where $1_{\Omega}$ is the indicator function of $\Omega\subseteq\mathbb{R}$.
\ By a limiting argument (from discrete to continuous), the PDF\ of $Z(t)$ as
$t\rightarrow\infty$ has Laplace transform \cite{Knd-que4}%
\[%
\begin{array}
[c]{ccc}%
\dfrac{(1-\rho)\theta}{\theta-\ln\left(  1+\rho\,\theta\right)  }, &  &
\operatorname{Re}(\theta)>0
\end{array}
\]
which Daniels \cite{Dnls-que4} inverted to yield%
\[%
\begin{array}
[c]{ccc}%
f(z)=-(1-\rho)%
%TCIMACRO{\dint \limits_{0}^{\infty}}%
%BeginExpansion
{\displaystyle\int\limits_{0}^{\infty}}
%EndExpansion
\,\dfrac{d}{dz}\dfrac{(z+w)^{w-1}\exp\left[  -(z+w)/\rho\right]  }{\rho
^{w}\Gamma(w)}\,dw, &  & z>0
\end{array}
\]
with a point mass $1-\rho$ at $z=0$. \ We seek an experimental approach to
verify this PDF. \ How might one efficiently simulate $Z(t)$ for suitably
large $t$?\ \ Offers of assistance would be most appreciated. \ We wonder too
if Prabhu's formula could possibly be reconfigured to play a role in this
inquiry. \ The fact that $v<\infty$ earlier but $v=\infty$ here is an issue;
the fact that $X_{t}$ was the precise inflow at time $t$ whereas $X(t)$ is an
accumulated inflow over $(0,t]$ is another issue. \ 

\section{Acknowledgements \ }

Khaled Hamed was so kind to answer several questions of mine; this paper would
not have been possible without his very helpful articles \cite{Hm1-que4,
Hm2-que4}. In particular, he appears to be the first author to specify the
role of the offset $\delta$ when $v$ is not an integer multiple of $m$. \ I am
grateful to innumerable software developers. \ The symbolic manipulations
described here are tailor-made for Mathematica, and the simulations employed
here to check predictions are ideal for R.

\end{document}